\documentclass[12pt]{article}
\usepackage{fullpage}
\usepackage{graphicx}
\usepackage{float}
\restylefloat{figure}

\def\RR{\hbox{I\kern-.2em\hbox{R}}}
\newcommand{\qed}{\hbox to 0pt{}\hfill$\rlap{$\sqcap$}\sqcup$
\vspace{3mm}}

\newtheorem{uess}{Lemma}
\newtheorem{guess}{Theorem}
\newtheorem{corollary}{Corollary}

\title{Nonoscillation  and Stability of the Second Order 
Ordinary Differential Equations with a Damping Term}
\author{Leonid Berezansky $^{1}$ \\Department of Mathematics,
Ben-Gurion University of the Negev, \\
Beer-Sheva 84105, Israel,
\\
Elena Braverman $^{2}$ \\
Department of Mathematics and Statistics, University of Calgary, \\
2500 University Drive N.W., Calgary, AB T2N 1N4, Canada\\
and   Alexander Domoshnitsky  $^{3}$ \\
Department of Mathematics and Computer Science,\\
Ariel University Center Samaria,
Ariel 44837, Israel}
\begin{document}
\maketitle

\begin{abstract}
In this paper we consider the linear
ordinary equation of the second order 
$$
\pounds x(t)\equiv 
\ddot{x}(t) +a(t)\dot{x}(t)+b(t)x(t)=f(t),  \eqno{(0.1)} 
$$
and the corresponding homogeneous equation 
$$
\ddot{x}(t) +a(t)\dot{x}(t)+b(t)x(t)=0.
\eqno{(0.2)} 
$$
Note that $[\alpha ,\beta ]$ is called a nonoscillation interval if every
nontrivial solution has at most one zero on this interval. Many
investigations which seem to have no connection such as
differential inequalities, the Polia-Mammana decomposition (i.e.
representation of the operator $\pounds$ in the form of products of the 
first order differential operators), unique solvability of the 
interpolation problems, kernels oscillation, separation of zeros, zones 
of Lyapunov's stability and some others have a certain common basis -
nonoscillation. Presumably Sturm was the first to consider the two
problems which naturally appear here: to develop corollaries of
nonoscillation and to find methods to check nonoscillation. In this paper we
obtain several tests for nonoscillation on the semiaxis and apply them to
propose new results on asymptotic properties and the exponential
stability of the second order equation (0.2).
Using the Floquet representations and upper and lower estimates of 
nonoscillation intervals of oscillatory solutions we deduce results on the 
exponential and Lyapunov's stability and instability of equation (0.2). 
\end{abstract}

{\bf Keywords and Phrases}: ordinary differential
equations, nonoscillation interval, exponential stability, boundary 
value problems, Cauchy function, Green's function, Floquet theory.

{\bf AMS(MOS) subject classification:} 34D, 34A30.

\thispagestyle{empty}

\section{Introduction}

This paper deals with the equation 
\begin{equation}
\label{1}
\ddot{x}(t) +a(t)\dot{x}(t)+b(t)x(t)=f(t),~ t\in[0,\omega],
\end{equation}
with locally summable coefficients $a,b,f,$ which together with nonlinear
equation 
\begin{equation}
\label{2}
\ddot{x}(t)=g(t,x(t),\dot{x}(t)),\;t\in [0,\omega ],
\end{equation}
continue to attract attention of many mathematicians due to their
significance in applications. In this paper we obtain two group of 
results: on nonoscillation and on exponential stability. To obtain  
stability results we will use nonoscillatory equations as the so-called 
model equations for the left or right regularization and then use the 
classical Bohl-Perron theorem \cite{DK,H}.
Another approach to stability study for periodic 
equations in the present paper is based on lower and upper estimates of 
the distance between two adjacent zeros (i.e., of nonoscillation 
intervals) for nontrivial solutions of homogeneous equations. The 
foundations of this approach can be found in the work by Zhukovskii 
\cite{zhukovskii}, Kre\u{\i}n \cite{krein} and Yakubovich
\cite{yakubovich1}. 
Zones of Lyapunov's stability can be also studied on this basis.
This explains why we connect the different 
areas together as well as the fact that actually our approach develops
applications of the classical nonoscillation area. Note that we obtain new
exponential stability conditions for equations with measurable 
coefficients. In most stability conditions it was assumed that $b(t)\equiv 
b>0$ \cite{H1,H2,LN,PS,S}, 
$b(t)\geq 0$ is a differentialble function \cite{BH,H3,I,KG,MV}
or some restrictions like slow varying coefficients \cite{DIS,G1,G2}
were imposed. We
consider here equation (\ref{1}) without usual restrictions on parameters of
the equations, the coefficients are even not required to be continuous.

Let us describe nonoscillation in general in order to understand how the
results of the present paper develop also many other topics. The 
nonoscillation area
consists of many topics which seem to have no relevance to each  
others but they are deeply connected.

The classical de la Vall\'{e}e-Poussin theorem claims that existence of a
positive function $v$ such that $v^{\prime \prime }(t)+p(t)v(t)\leq 0$
for $t\in [ 0,\omega ]$ impiles that $[0,\omega ]$ is a
nonoscillation interval. The idea of theorems on differential 
inequalities can be formulated as follows: under certain 
conditions solutions of inequalities are 
greater or less than the solution of the equation. The idea to construct 
an approximate
integration method for the numerical solution of differential equations  
based on the comparison of solutions of equations and inequalities first 
appeared in the works of famous Russian mathematician 
Chaplygin~\cite{chaplygin} 
and later was developed by other famous Russian mathematician Luzin~\cite
{luzin}. Concerning our object we can formulate the differential 
inequality theorem in the form: {\it under certain conditions the 
inequalities} 
\begin{equation}
\label{3}
(\pounds y)(t)\geq (\pounds x)(t),~t\in [ 0,\omega ],\,\;y(0)\geq
x(0),\;y^{\prime }(0)\geq x^{\prime }(0)\;
\end{equation}
{\it imply }$y(t)\geq x(t)${\it \ for }$t\in [0,\omega ].$

Independly Azbelev~\cite{azbelev1}, Beckenbach, 
Bellman~\cite{bellman} and Wilkins~\cite{wilkins} established that 
(\ref{3}) can be applied in a nonoscillation interval $[0,\omega ]$ only. 
The general solution of equation (\ref{1}) has the following 
representation 
\begin{equation}
\label{4}
x(t)=\int_{0}^{t}C(t,s)f(s)ds+x_{1}(t)x_{0}+x_{2}(t)x_{0}^{\prime },
\end{equation}
where 
$\ x_{1}$ and $x_{2}$ are the solutions of the homogeneous equation 
\begin{equation}
\label{5}
\pounds x(t)\equiv \ddot{x}(t)+a(t)\dot{x}(t)+b(t)x(t)=0,\;t\in [0,\omega ],
\end{equation}
satisfying the initial conditions $x_{1}(0)=1$, $x_{1}^{\prime }(0)=0$ and 
$x_{2}(0)=0$, $x_{2}^{\prime }(0)=1$, respectively,
$x_{0}$ and$\ x_{0}^{\prime }$ are corresponding constants.
The kernel $C(t,s)$ of the integral in solution's representation (\ref{4}) 
is called {\em the Cauchy function} of equation (\ref{1}). {\em The 
fundamental function} $X(t,s)$ of (\ref{1}) is defined as follows: 
$X(t,s)$ for each fixed $s\geq 0$ as a function of $t$ satisfies 
\begin{equation}
\label{6}
\ddot{x}(t)+a(t)\dot{x}(t)+b(t)x(t)=0,\;t\in [ s,\omega],
\end{equation}
\begin{equation}
\label{7}
x(s)=0,\;x^{\prime }(s)=1.
\end{equation}
For equation (\ref{1}) the Cauchy function and the fundamental function 
$X(t,s)$ coincide  \cite{AMR}.
We assume that $X(t,s)=C(t,s)=0$ for $0\leq t<s.$ 

If the solution of the initial value problem $
(\pounds x)(t)=0,\;x(0)=0,\;x^{\prime }(0)=1$ does not vanish at the point $%
t=\omega $\ , then the boundary value problem 
\begin{equation}
\label{8}
\pounds x(t)\equiv \ddot{x}(t)+a(t)\dot{x}(t)+b(t)x(t)=f(t),\;t\in [ 0,\omega ],
x(0)=0,\;x(\omega )=0,%
\end{equation}
is uniquely solvable and its solution has the representation 
\begin{equation}
\label{9}
x(t)=\int_{0}^{\omega }G(t,s)f(s)ds,
\end{equation}
where the kernel of the integral representation $G(t,s)$ is called the
Green's function of the problem (\ref{8}). Differential inequality 
theorems are actually results on positivity or negativity of corresponding 
Green's functions.

The equivalence of nonoscillation and the unique solvability of the
interpolation problems 
\begin{equation}
\label{10}
\begin{array}{l}{\displaystyle
\pounds x(t)\equiv \ddot{x}(t)+a(t)\dot{x}(t)+b(t)x(t)=f(t),~~t\in 
[0,\omega ],}\\
x(t_{1})=0,~~
\;x(t_{2})=0, ~~~0\leq t_{1}<t_{2}\leq \omega ,
\end{array}
\end{equation}
is obvious \cite{levin}. Let us say that the Green's function of problem 
(\ref{10}) behaves regularly if 
\begin{equation}
\label{11}
G(t,s)(t-t_{1})(t-t_{2})>0,\;t,s\in (0,\omega ),\;t\neq t_{1},\;t\neq t_{2}.%
\end{equation}
It was first proven in \cite{chichkin} that nonoscillation of the equation $
\pounds x=0$\ on the interval $[0,\omega ]$\ is necessary and sufficient for
the regular behavior of  Green's functions of interpolation problem 
(\ref{10})
(see also the well known paper by Levin \cite{levin}). 

Let us denote by $C_{[0,\omega ]}$ the space of continuous functions $
x:[0,\omega ]\rightarrow \RR$ with the norm $\left\| x\right\| =\max_{t\in
[0,\omega ]}\left| x(t)\right| ,$ and by $D_{[0,\omega ]}$ the linear
space of functions $x:[0,\omega ]\rightarrow \RR$ with absolute 
continuous  
derivatives. Let $b(t)=b^{+}(t)-b^{-}(t),$ where $b^{+}(t)\geq
0,\;b^{-}(t)\geq 0,$ and 
\begin{equation}
\label{12}
\pounds _{0}x(t)\equiv \ddot{x}(t)+a(t)\dot{x}(t)-b^{-}(t)x(t)=0,\;t\in [ 
0,\omega ].
\end{equation}

It was proven in \cite{AzbelevDom2} that the boundary value problem 
\begin{equation}
\label{13}
\pounds _{0}x(t)\equiv \ddot{x}(t)+a(t)\dot{x}(t) - b^{-}(t)x(t)=f(t),
~~t\in [ 0,\omega ],
x(0)=0, ~x(\omega)=0,
\end{equation}
is uniquely solvable and 
its Green's function $G_{0}(t,s)$ is negative 
$G_{0}(t,s)<0$ for $t,s\in (0,\omega )$. In the
space $C_{[0,\omega ]}$ let us define the integral operator 
$K:C_{[0,\omega]}\rightarrow C_{[0,\omega ]}$ by the equality 
\begin{equation}
\label{14}
(Kx)(t)=-\int_{0}^{\omega }G_{0}(t,s)b^{+}(s)x(s)ds.
\end{equation}
Note that the operator $K$ actually maps $C_{[0,\omega]}$ into $
D_{[0,\omega ]}$ due to the properties of Green's function $G_{0}(t,s).$
Hence the equation $x=Kx+g$, where $g(t)=\int_{0}^{\omega }G_{0}(t,s)f(s)ds$, 
is equivalent to the boundary value problem (\ref{8}).

The above argument can be summarized in the form of the statement 
on six equivalences.

{\bf Theorem A} \cite{AzbelevDom2}. {\it The following assertions are
equivalent:}

{\it 1) there is }$v\in D_{[0,\omega ]}$     
{\it \ such that }$v(t)\geq 0${\it \ and }$(
\pounds v)(t)\leq 0${\it \ for }$t\in [ 0,\omega ]${\it \ and }$%
v(0)+v(\omega )-\int_{0}^{\omega }(\pounds v)(t)dt>0;$

{\it 2) }$C(t,s)>0${\it \ for }$0\leq s\leq t\leq \omega ;$

{\it 3) problem (\ref{8}) has a unique solution for each summable }
$f${\it \ and }$G(t,s)<0${\it \ for }$t,s\in (0,\omega );$

{\it 4) each nontrivial solution of the homogeneous equation \pounds }$x=0$%
{\it \ has at most one zero on }$[0,\omega ];$

{\it 5) the spectral radius of the operator }$K:C_{[0,\omega ]}\rightarrow
C_{[0,\omega ]}${\it \ is less than one;}

{\it 6) problem (\ref{10}) has a unique solution for each summable }$f${\it \ and
its Green's function behaves regularly for }$t,s\in (0,\omega ).$

\bigskip 

The Polia-Mammana decomposition \cite{Mammana, Polia} is a possibility for
representation of the operator $\pounds :D_{[0,\omega ]}\rightarrow
L_{[0,\omega ]}$\ in the form of products of the first order differential
operators 
\begin{equation}
\label{15}
\pounds =h_{2}\frac{d}{dt}h_{1}\frac{d}{dt}h_{0},
\end{equation}
where the real valued functions $h_{i}$ do not have zeros and are smooth
enough. This decomposition is possible if and only if $[0,\omega ]$ is a
nonoscillation interval. In particular, representation (\ref{15}) allows 
us to obtain the generalized Rolle's theorem \cite{levin}: {\em if the 
solution $x$ of equation (\ref{1}) has more than 3 zeros on the 
nonoscillation interval, then the function $f$ has at least one 
zero on this interval}.

The theory of oscillatory kernels plays an important role in oscillation of
mechanical systems \cite{Gantmakher}. The oscillatory kernel $G(t,s)$
is characterized by the inequalities 
\begin{equation}
\label{16}
G(t,s)>0,\;\det \left| G(t_{i},s_{j})\right| _{1}^{m}\geq
0,\;0<t_{1}<...<t_{m}<\omega ,\;0<s_{1}<...<s_{m}<\omega ,\;m=1,2,...
\end{equation}
while for $t_{i}=s_{i}$\ $(1=1,2,,...,m)$ the inequality 
in (\ref{16}) has to be strict. 
In \cite{Gantmakher} it was demonstrated that actually the fact 
that the kernel is oscillatory means that the integral operator 
$G:L_{[0,\omega]}\rightarrow D_{[0,\omega ]}$ of the form 
$(Gf)(t)=\int_{0}^{\omega
}G(t,s)f(s)ds$\ does not increase the number of sign's changes. If we
consider the integral operator, where the kernel $G(t,s)$
is the Green's function of problem (\ref{8}), then in this case the inverse
operator $\pounds :D_{[0,\omega ]}\rightarrow L_{[0,\omega ]}$\ defined on
the functions satisfying the conditions $x(0)=0,\;x(\omega )=0$ should not
decrease the number of sign's changes. Although a direct verification of 
an infinite number of inequalities (\ref{16}) is possible only for very 
simple kernels and cannot be implemented in most interesting for 
applications cases, non-decreasing of the number of 
sign's changes for the integral operators with Green's functions as  
kernels can be checked
through the Polia-Mammana decomposition and the generalized Rolle's theorem.
This connection of oscillatory kernels and Polia-Mammana decomposition was
discovered by M.G.Kre\u{\i}n. 

Thus nonoscillation solves in many important
cases the problem of checking oscillatory kernels.
As a conclusion we note that each new nonoscillation result or test 
develops all these directions.

\section{Preliminaries}
\setcounter{equation}{0}

Let us start with the following simple corollaries of Theorem A.
Choosing $v(t)=\exp (\lambda t)$ in the first assertion of Theorem A, we 
obtain the following result.

\begin{corollary}
 If there exists such a real constant $\lambda $ such that 
\begin{equation}
\label{19}
\lambda ^{2}+a(t)\lambda +b(t)\leq 0,
\end{equation}
 then for each $\omega $  assertion 2)-6) of Theorem A are true.
If in addition $b(t)\geq \beta >0$ and this $\lambda < 0,$  then
equation (\ref{1}) is  exponentially stable (the fundamental function has an
exponential estimate).
\end{corollary}
\bigskip

Solving the equation 
\begin{equation}
\label{20}
\lambda ^{2}+a(t)\lambda +b(t)=0,
\end{equation}
for all $t,$ we get $\lambda _{1}(t)=-\frac{a(t)}{2}-\sqrt{\frac{a^{2}(t)}{4}%
-b(t)}$ and $\lambda _{2}(t)=-\frac{a(t)}{2}+\sqrt{\frac{a^{2}(t)}{4}-b(t)}.$

\begin{corollary}$ \cite{levin}$
 If for sufficiently large $t\geq
t_{0}$ there exist $\nu _{0},\nu _{1}$ and $\nu _{2}$ such
that $\lambda _{1}(t)$  and $\lambda _{2}(t)$  are real functions
and $\nu _{0}\leq \lambda _{1}(t)\leq \nu _{1}\leq \lambda
_{2}(t)\leq \nu _{2},$ where $\nu _{0}<\nu _{1}<\nu _{2},$ then
the fundamental system of equation (\ref{1}) satisfies the inequalities 
\begin{equation}
\label{21}
c_{i}\exp (v_{i-1}t)\leq x_{i}(t)\leq d_{i}\exp
(v_{i}t),\;i=1,2\;(c_{i},d_{i}>0;\;t\geq t_{0}).
\end{equation}
If in addition $\nu_2 <0,$  then equation (\ref{1}) is  exponentially stable.
\end{corollary}

In order to  study stability properties, we consider the scalar  
differential equation of the  second order (\ref{1}) under the following 
conditions:

(a1) $a(t), b(t)$ are Lebesgue measurable and  
essentially bounded functions on $[0,\infty)$.

(a2) $f:[t_0,\infty)\rightarrow R $ is a Lebesgue measurable locally 
essentially bounded function. 

\noindent

\noindent
{\bf Definition}. Eq.~(\ref{1}) is {\em (uniformly) exponentially stable}, 
if the fundamental function  $X(t,s)$ of (\ref{1}) {\em has
an exponential estimate}
if there exist
positive numbers 
$K>0, \lambda>0$, such that
\begin{equation}
\label{17}
|X(t,s)|\leq K~e^{-\lambda (t-s)},~~t\geq s\geq 0.
\end{equation}

Consider the equation
\begin{equation}
\label{18}
\ddot{x}(t) +a\dot{x}(t)+bx(t)=0, 
\end{equation}
where $a>0, b>0$ are positive numbers. 
This equation is exponentially stable.
Denote by $Y(t,s)$ the fundamental function of  (\ref{18}).

\begin{uess}
Let $a>0, b>0$.

1) If  $a^2>4b$ then ${\displaystyle \int_0^t|Y(t,s)|ds\leq\frac{1}{b}, ~~
\int_0^t|Y^{'}_t(t,s)|ds\leq\frac{2a}{\sqrt{a^2-4b}(a-\sqrt{a^2-4b})} }$.
\vspace{2mm}

2) If $a^2<4b$ then ${\displaystyle \int_0^t|Y(t,s)|ds\leq 
\frac{4}{a\sqrt{4b-a^2}}, ~~
\int_0^t|Y^{'}_t(t,s)|ds\leq \frac{2(a+\sqrt{4b-a^2})}{a\sqrt{4b-a^2}} }$.
\vspace{2mm}

3)  If  $a^2=4b$ then ${\displaystyle \int_0^t|Y(t,s)|ds\leq \frac{1}{b}, 
~~
\int_0^t|Y^{'}_t(t,s)|ds\leq \frac{2}{\sqrt{b}}}$.
\end{uess}
{\bf Proof.}
If $a^2>4b$ then the characteristic equation $\lambda^2+a\lambda+b=0$
has two negative roots 
$$
\lambda_{1,2}=\frac{-a\pm \sqrt{a^2-4b}}{2}, ~0\geq \lambda_1>\lambda_2.
$$
By simple calculation we have
$$
0<Y(t,s)=\frac{1}{\sqrt{a^2-4b}}\left(e^{\lambda_1(t-s)}-e^{\lambda_2(t-s)}
\right).
$$
Since ${\displaystyle 
\frac{1}{\lambda_1}e^{\lambda_1 t} < \frac{1}{\lambda_2}e^{\lambda_2 t} }$
then 
$$
\int_0^t Y(t,s)ds\leq \frac{1}{\sqrt{a^2-4b}}
\left(\frac{1}{\lambda_2}-\frac{1}{\lambda_1}\right)=\frac{1}{b}.
$$
We have $Y_t^{'}(t,s)=\frac{1}{\sqrt{a^2-4b}}\left(\lambda_1e^{\lambda_1(t-s)}-\lambda_2e^{\lambda_2(t-s)}\right)$.
Hence 
${\displaystyle |Y_t^{'}(t,s)|\leq 
\frac{|\lambda_1|+|\lambda_2|}{\sqrt{a^2-4b}}e^{\lambda_1(t-s)} }$
$$=\frac{a}{\sqrt{a^2-4b}}e^{\lambda_1(t-s)} \mbox{~~and~~~}
\int_0^t|Y^{'}_t(t,s)|ds\leq \frac{2a}{\sqrt{a^2-4b}(a-\sqrt{a^2-4b})}.$$

If $a^2-4b<0$, then the characteristic equation has two complex roots 
and the fundamental function has the form 
$$
Y(t,s)=\frac{2}{\sqrt{4b-a^2}}e^{-\frac{a}{2}(t-s)}
\sin\left( \frac{\sqrt{4b-a^2}}{2}(t-s) \right).
$$
Hence
${\displaystyle
|Y(t,s)|\leq \frac{2}{\sqrt{4b-a^2}}e^{-\frac{a}{2}(t-s)}
}$
and
${\displaystyle
\int_0^t|Y(t,s)|ds\leq \frac{4}{a\sqrt{4b-a^2}}.
}$
The second inequality in 2) is proven in a similar way as the second 
inequality in 1).

If $a^2=4b$, then the characteristic equation has the double root 
$\lambda=-\frac{a}{2}$  and 
$$
Y(t,s)=(t-s)e^{-\frac{a}{2}(t-s)}.
$$
Since
${\displaystyle \int_0^{\infty} s e^{-as/2}= \frac{4}{a^2}=\frac{1}{b} }$
then the first inequality in 3) holds. The second inequality in 3) is 
proven similarly to the previous cases.
\qed

Let us introduce some functional spaces on a semi-axis. Denote
by $L_{\infty}[t_0,\infty)$ the space of all essentially bounded on
$[t_0,\infty)$
functions and by $C[t_0,\infty)$ the space of all continuous
bounded on  $[t_0,\infty)$ functions
with the supremum norm. 

\begin{uess}$\cite{DK,H}$
Suppose  (a1)-(a2) hold and there exists $t_0\geq 0$ 
such that for every 
$f\in L_{\infty}[t_0,\infty)$ 
both the solution $x$ of the problem 
$$
\ddot{x}(t) +a(t)\dot{x}(t)+b(t)x(t)=f(t),~ t\geq t_0,
$$$$
x(t)=0, \dot{x}(t)=0,~t\leq t_0,
$$
and its derivative $\dot{x}$ belong to $C[t_0,\infty)$. Then 
equation (\ref{1}) is exponentially
stable.
\end{uess}

\section{Integro-differential equation}
\setcounter{equation}{0}

To obtain positiveness conditions for the fundamental function of equation 
(\ref{1}) we consider first a similar problem for the following 
integro-differential equation
\begin{equation}
\label{22}
\dot{y}(t)+\int_0^t e^{-\int _s^ta(\xi)d \xi}b(s)y(s)ds=0,
\end{equation}
for which we assume that condition (a1) holds.

Together with (\ref{8}) we consider 
for each $t_0\geq 0$ the initial value problem
\begin{equation}
\label{23}
\dot{y}(t)+\int_{t_0}^t e^{-\int _s^ta(\xi)d \xi}b(s)y(s)ds=f(t), 
\end{equation}
\begin{equation}
\label{24}
y(t_0)=y_0.
\end{equation}
We assume that condition (a2) holds for the function $f(t)$.

\noindent
{\bf Definition.} A locally absolutely
continuous on $[t_0,\infty)$ function $y:\RR \rightarrow \RR$ is called 
{\em a solution} of  problem (\ref{23}), (\ref{24}) 
if it satisfies equation (\ref{23}) for almost every 
$t\in [t_0, \infty)$, and equality (\ref{24}) for $t= t_0$.

\noindent
{\bf Definition.} For each $s\geq 0$ the solution $Y(t,s)$ of the problem
\begin{equation}
\label{25}
\dot{y}(t)+\int_{s}^t e^{-\int _s^ta(\xi)d \xi}b(\tau)y(\tau)d\tau=0, ~~
t>s, ~~ y(s)=1,
\end{equation}
is called {\em a fundamental function} of equation (\ref{22}).

We assume $Y(t,s)=0, t<s$.

\begin{uess}$\cite{H}$
Let (a1)-(a2) hold. Then there exists one and only one
solution of  problem (\ref{23}), (\ref{24}) that can be presented in the form
\begin{equation}
\label{26}
y(t)=Y(t,t_0)y_0+\int_{t_0}^t Y(t,s)f(s)ds.
\end{equation}
\end{uess}

Let us obtain conditions under which equation (\ref{22}) has a 
positive solution.
We remark that
the theorem remains true if the zero initial point is replaced by any 
$t_0 \geq 0$.

In future we will apply the following result (\cite{AMR}, Theorem 1.2.1).
Let $L[0,T]$ be the space of all integrable on $[0,T]$ functions with the 
norm $\|f\|=\int^{T}_{0}|f(s)|ds$.

\begin{uess}$\cite{AMR}$
\label{lemmaAzb}
Suppose that a function $r(t,s)$ is measurable over the square 
$[a,b]\times[a,b]$, for almost every $s$
the function $r(t,s)$ as a function of $t$ has 
finite one-sided limits 
at each point $t$ 
and there exists a function
$v\in L[0,T]$ such that $|r(\cdot,s)|\leq v(\cdot)$. Then the integral 
operator 
$$
(Ky)(t)=\int_a^b r(t,s)y(s)ds
$$
maps  L[0,T] to L[0,T] and is a compact operator in this space.
\end{uess}
\begin{guess}
Suppose $a(t)\geq 0, b(t)\geq 0$. Then
the following conditions are equivalent for equation (\ref{22}).

1) There exists a positive solution of the inequality
\begin{equation}
\label{27}
\dot{y}(t)+\int_{0}^t e^{-\int_s^t a(\xi)d \xi}b(s)y(s)ds\leq 0, ~t \geq 0.
\end{equation}

2) There exists a locally essentially bounded  function $u(t)\geq 0$ such that
\begin{equation}
\label{28}
u(t)\geq \int_0^t e^{-\int_s^t [a(\xi)-u(\xi)]d \xi}b(s)ds, ~t \geq 0.
\end{equation}

3) $Y(t,s)>0, t\geq s\geq 0$.

4) There exists a positive solution of  equation (\ref{8}) for $t \geq 0$.
\end{guess}
{\bf Proof.}
1) $\Longrightarrow$ 2). 
Suppose $y(t)>0$ is a solution of (\ref{27}). 
Hence on any bounded interval $y(t)\geq \alpha>0$, then
${\displaystyle u(t)=-\frac{\dot{y}(t)}{y(t)}\geq 0}$  and is an 
essentially locally bounded function. We also have 
$$
y(t)=y(0)e^{-\int_0^tu(\xi)d \xi}, ~~
\dot{y}(t)=-y(0)u(t)e^{-\int_0^tu(\xi)d \xi}.
$$
Substituting $y, \dot{y}$ in (\ref{27}) we have
\begin{eqnarray*}
& & \dot{y}(t)+  \int_{0}^t e^{-\int_s^t a(\xi)d \xi}b(s)y(s)ds
\\
& = & -y(0)u(t)e^{-\int_0^t u(\xi)d \xi}+y(0)\int_0^t 
e^{-\int_s^t a(\xi)d \xi}b(s)e^{-\int_0^su(\xi)d \xi}ds
\\
& = & -y(0)e^{-\int_0^tu(\xi))d \xi}\left[u(t)-
\int_0^t e^{-\int_s^t [a(\xi)-u(\xi)]d \xi}b(s)ds\right]\leq 0.
\end{eqnarray*}
Hence $u(t)\geq 0$ is a solution of inequality (\ref{28}).

2) $\Longrightarrow$ 3).
Consider the nonhomogeneous equation (\ref{27})
\begin{equation}
\label{29}
\dot{y}(t)+\int_{0}^t e^{-\int_s^t a(\xi)d \xi}b(s)y(s)ds=f(t),~~ y(0)=0.
\end{equation}
A solution of  (\ref{29})  will 
be sought in the form
\begin{equation}
\label{30}
y(t)=\int_{0}^t e^{-\int_s^t u(\xi)d \xi}z(s)ds,~~
\dot{y}(t)=z(t)-u(t)\int_{0}^t e^{-\int_s^t u(\xi)d \xi}z(s)ds,
\end{equation}
so $\dot{y}(t)+u(t)y(t)=z(t), y(0)=0.$

After substituting (\ref{30}) into (\ref{29}) we have
\begin{equation}
\label{31}
z(t)-u(t)\int_{0}^t e^{-\int_s^t u(\xi)d \xi}z(s)ds
+\int_{0}^t e^{-\int_s^t a(\xi)d \xi}b(s)
\int_{0}^s e^{-\int_{\tau}^s u(\xi)d \xi}z(\tau)d\tau ds=f(t).
\end{equation}
Since 
\begin{eqnarray*}
& & \int_{0}^t \left[ e^{-\int_s^t a(\xi)d \xi}b(s)\int_{0}^s 
e^{-\int_{\tau}^s
u(\xi)d \xi}z(\tau)d\tau \right] ds
\\
& = & \int_{0}^t \left[ \int_{\tau}^t e^{-\int_s^t a(\xi)d \xi}b(s)
e^{-\int_{\tau}^su(\xi)d \xi}ds \right] z(\tau)d\tau
\\
& = & \int_{0}^t e^{-\int_{\tau}^t u(\xi)d \xi} \left[ \int_{\tau}^t
e^{-\int_s^t (a(\xi)-u(\xi))d \xi }b(s)ds \right] z(\tau)d\tau,
\end{eqnarray*}
then equation (\ref{31}) can be rewritten as
\begin{equation}
\label{32}
z(t)-\int_{0}^t e^{-\int_{\tau}^t u(\xi)d \xi}
\left[u(t)-\int_{\tau}^t e^{-\int_s^t (a(\xi)-u(\xi))d \xi }b(s)ds\right]
z(\tau)d\tau=f(t).
\end{equation}
On every finite interval $[0,T]$ equation (\ref{32}) has the form
\begin{equation}
\label{33}
z-Hz=f, ~t\in [0,T].
\end{equation}
Operator $H:L[0,T]\rightarrow L[0,T]$ is bounded.
In order to show that this operator is compact we apply  
Lemma \ref{lemmaAzb}.
Operator $H$ can be rewritten in the form $H=PH_1-H_2$, where
$$
(Pz)(t)=u(t)z(t),~(H_1z)(t)=\int_{0}^t e^{-\int_{\tau}^t 
a(\xi)d \xi}z(\tau)d\tau,
$$$$
(H_2z)(t)=\int_{0}^t 
\left[ e^{-\int_{\tau}^t a(\xi)d \xi}\int_{\tau}^t e^{-\int_s^t 
(a(\xi)-u(\xi))d \xi }b(s)ds \right] z(\tau)d\tau.
$$
It is easy to see that for operators $H_1, H_2$ all conditions of Lemma 5 hold.
Then these operators are compact.
Operator $P$ is bounded operator, 
hence operator $H$ is a compact Volterra integral operator
with spectral radius $r(T)=0$ \cite{R}.
Hence for the solution of equation (\ref{33}) we have
$z=(I-H)^{-1}f$, where $I$ is the identity operator.

Since 
$$
u(t)-\int_{\tau}^t e^{-\int_s^t (a(\xi)-u(\xi))d \xi }b(s)ds
\geq u(t)-\int_0^t e^{-\int_s^t (a(\xi)-u(\xi))d \xi }b(s)ds\geq 0,
$$
then $H$ is a positive operator.
Hence $(I-H)^{-1}= 1+H+H^2+H^3+ \dots $
is also a positive operator.

Suppose now that in the equation (\ref{33}) we have $f(t)\geq 0$. Then for 
the solution of (\ref{33}) 
we have $z(t)\geq 0$. Equality (\ref{16}) implies that for every 
right-hand side $f(t)\geq 0$ the solution of equation (\ref{29}) $y(t)\geq 0$. 
But $y(t)=\int_0^t Y(t,s) f(s)ds$. Hence $Y(t,s)\geq 0,~ 0\leq s\leq t\leq T$.
Since $T>0$ is an arbitrary number then $Y(t,s)\geq 0,~ 0\leq s\leq t<\infty$. 
We only have to prove that the strong inequality for $Y(t,s)>0$ holds.

After substituting $y(t)=e^{-\int_0^tu(\xi)~d \xi} $ in the left-hand 
side of 
equation (\ref{29}) we see that this function is a solution of (\ref{29})
with $f(t)<0$. By the solution representation formula (\ref{26}) we have
$$
y(t)=Y(t,0)+\int_0^t Y(t,s) f(s)ds.
$$
Hence  $Y(t,0)\geq y(t)>0$. The general case $Y(t,s)>0$ is considered similarly.

Implications 3)$\Longrightarrow$4) and 4)$\Longrightarrow$1) are evident. 
\qed

\noindent
{\bf Remark.} Nonoscillation conditions for general integro-differential 
equations with a bounded memory were obtained in \cite{BB}.
However, these results are not applicable to equation (\ref{22}).
Nonoscillation results for integro-differential 
equation can also be found in \cite{DG2001}.

\section{Positive Solutions}
\setcounter{equation}{0}

The following lemma gives a connection between 
equations (\ref{1}) and (\ref{22}).

\begin{uess}
Denote by $x_1(t), x_2(t)$ and $X(t,s)$
the fundamental system and the fundamental function of (\ref{1}), 
respectively, by $Y(t,s)$
the fundamental function of (\ref{8}). Then 
$$
x_1(t)=Y(t,0), ~x_2(t)=\int_0^t Y(t,\tau)e^{-\int_0^{\tau}a(\xi)d\xi}d\tau,
~X(t,s)=\int_s^t Y(t,\tau)e^{-\int_s^{\tau}a(\xi)d\xi}d\tau.
$$ 
\end{uess}
{\bf Proof.}
For the solution of equation (\ref{1}) we have
$$
\dot{x}(t)=e^{-\int_0^{t}a(\xi)d\xi}\dot{x}(0)
-\int_0^{t}e^{-\int_s^{t}a(\xi)d\xi}b(s)x(s)ds.
$$
Hence
$$
\dot{x}(t)+\int_0^{t}e^{-\int_s^{t}a(\xi)d\xi}b(s)x(s)ds
=e^{-\int_0^{t}a(\xi)d\xi}\dot{x}(0)
$$
and
$$
x(t)=Y(t,0)x(0)+\left[ \int_0^t 
Y(t,\tau)e^{-\int_0^{\tau}a(\xi)d\xi}d\tau\right] \dot{x}(0).
$$
But for the solution of (\ref{1}) we have another representation:
$$
x(t)=x_1(t)x(0)+x_2(t)\dot{x}(0).
$$
The equalities for the fundamental system of (\ref{1}) are proven. 
Since $X(t,s)=x_2(t,s)$, then the proof of the equality for $X(t,s)$ is 
similar.
\qed

\begin{corollary}
If the fundamental function $Y(t,s)$ of (\ref{22}) is positive then
the fundamental system and the fundamental function of (\ref{1}) are positive.  
\end{corollary}
\begin{corollary}
Suppose $a(t)\geq 0, b(t)\geq 0,$ and the fundamental function of (\ref{8}) is positive.
Then 
\begin{equation}
\label{34}
0\leq \int_{t_0}^t X(t,s)b(s)ds\leq 1,
\end{equation}
where $X(t,s)$ is the fundamental function of (\ref{1}).
\end{corollary}
{\bf Proof.}
The function $x(t)\equiv 1$ is the solution of (\ref{2}) with $f(t)=b(t)$.
By the solution representation formula we have
$$
1=x_1(t)+\int_{t_0}^t X(t,s)b(s)ds.
$$
Corollary 3 implies $x_1(t)>0, X(t,s)>0.$ Hence the inequality (\ref{34}) 
is valid.
\qed

Together with (\ref{1}) consider the following equation
\begin{equation}
\label{35}
\ddot{x}(t) +a_1(t)\dot{x}(t)+b_1(t)x(t)=0,
\end{equation}
where for $a_1(t), b_1(t)$ condition (a1) holds.
\begin{corollary}
Suppose $a_1(t)\geq a(t)\geq 0, b(t)\geq b_1(t)\geq 0,$
and equation  (\ref{22}) has a positive solution.
Then the fundamental function and the fundamental system of  (\ref{35})
 are positive.
\end{corollary}
{\bf Proof.}
If (\ref{22}) has a positive solution, then inequality (\ref{28}) has 
a nonnegative solution $u(t)\geq 0$. This function is a nonnegative solution of
inequality (\ref{28}) where $a(t)$ and $b(t)$ are replaced by $a_1(t)$ and 
$b_1(t)$.
Corollary 3 implies this corollary.
\qed

\begin{guess}
Suppose $a_1(t)\geq a(t)\geq 0, b(t)\geq b_1(t), b(t)\geq 0$.
If the fundamental function of (\ref{22}) is positive, 
then the fundamental function and the fundamental system of (\ref{35}) 
are positive.
\end{guess}
{\bf Proof.}
Consider the equation
\begin{equation}
\label{36}
\dot{y}(t)+\int_{0}^t e^{-\int_s^t a_1(\xi)d \xi}b_1(s)y(s)ds=f(t),~~ 
y(0)=0,
\end{equation}
which can be rewritten in the form 
\begin{eqnarray}
& & \dot{y}(t)+\int_{0}^t e^{-\int_s^t a(\xi)d \xi}b(s)y(s)ds
\nonumber
\\
& = &\int_{0}^t \left( e^{-\int_s^t a(\xi)d \xi}-e^{-\int_s^t
a_1(\xi)d \xi}\right)b(s)y(s)ds
\nonumber
\\
\label{37}
& + & \int_{0}^t e^{-\int_s^t a_1(\xi)d \xi}(b(s)-b_1(s))y(s)ds+f(t), 
 ~~ y(0)=0.
\end{eqnarray}
Hence (\ref{37})  is equivalent to the equation
\begin{equation}
\label{38}
y(t)=\int_{0}^t Y(t,s)\int_{0}^s \left(e^{-\int_{\tau}^s
a(\xi)d \xi}-e^{-\int_{\tau}^s
a_1(\xi)d \xi}\right)b(\tau)y(\tau)d\tau ~ ds
\end{equation}
$$
+\int_{0}^t Y(t,s)\int_{0}^s e^{-\int_{\tau}^s
a_1(\xi)d \xi}(b(\tau)-b_1(\tau))y(\tau)d\tau ~ ds + g(t),
$$
where $g(t)=\int_{0}^t Y(t,s) f(s)ds\geq 0$. 

Equation (\ref{38}) has the form $y=Hy+g$, where $H$ is a positive compact 
Volterra operator in the space $L{[0,T]}$
for every $T>0$. 
Hence the spectral radius $r(H)=0$ 
in the space $L{[0,T]}$ and then for the solution of (\ref{38})
we have $y(t)=((I-T)^{-1}g)(t)\geq 0$. The solution of 
the equation (\ref{36}) has the form
$y(t)=\int_{0}^t Y_1(t,s)f(s)ds$, where $Y_1(t,s)$ is the 
fundamental function of the equation
$$
\dot{y}(t)+\int_{0}^t e^{-\int_s^t a_1(\xi)d \xi}b_1(s)y(s)ds=0.
$$
We obtained that for every $f(t)\geq 0$ the solution of (\ref{36}) 
is also nonnegative.
Hence the inequality $Y_1(t,s)\geq 0$ holds.
Now, let us prove the strict inequality $Y_1(t,s)> 0$.

The function $Y_1(t,0)$ is the solution of the equation
\begin{equation}
\label{39}
\dot{y}(t)+\int_0^t e^{-\int_s^t a_1(\xi)d \xi}b_1(s)y(s)ds=0,~~ y(0)=1.
\end{equation}
After rewriting (\ref{39}) in the form
\begin{eqnarray}
& & \dot{y}(t) +\int_{0}^t e^{-\int_s^t a(\xi)d \xi}b(s)y(s)ds
\nonumber
\\
& = & \int_{0}^t \left(e^{-\int_s^t a(\xi)d \xi}-e^{-\int_s^t
a_1(\xi)d \xi}\right)b(s)y(s)ds
\nonumber
\\
\label{40}
& + & \int_{0}^t e^{-\int_s^t a_1(\xi)d \xi}(b(s)-b_1(s))y(s)ds,~~ y(0)=1,
\end{eqnarray}
equation (\ref{40}) has the form
$$
\dot{y}(t)+\int_{0}^t e^{-\int_s^t a(\xi)d \xi}b(s)y(s)ds =f(t), y(0)=1,
$$
where $f(t)\geq 0$. Hence
$$
Y_1(t,0)=Y(t,0)+\int_{0}^t Y(t,s)f(s)ds.
$$
Then $Y_1(t,0)\geq Y(t,0)>0$.

The general case $Y_1(t,s)>0$ is considered similarly.
Now, Lemma 6 implies the statement of the theorem.
\qed

\noindent
{\bf Remark.} The theorem remains true if we replace the zero initial 
point $0$ by any $t_0>0$.

Denote $b^+=\max\{b,0\}$.
\begin{corollary}
Suppose $a(t)\geq 0$ and the fundamental 
function of the equation
$$
\dot{y}(t)+\int_{0}^t e^{-\int_s^t a(\xi)d \xi}b^+(s)y(s)ds=0
$$
is positive. Then the fundamental function and the fundamental system of equation (\ref{1})
are positive.
\end{corollary}

Now let us proceed to explicit sufficient conditions when 
equation (\ref{1}) has a positive solution.

\begin{guess}
Suppose  at least one of the following conditions holds
\vspace{2mm}

1) $ a(t)\geq \int_{t_0}^t b^+(s)ds$,
\vspace{2mm}

2) there exists $\lambda>0$ such that  $a(t)\geq \lambda
b^+(t)+\frac{1}{\lambda},~ t\geq t_0$,  
\vspace{2mm}

3) $a(t)\geq 0$ and there exists $\lambda>0$ such that $\int_{t_0}^t e^{-\int_s^t
(a(\xi)-\lambda)d\xi}b^+(s)ds\leq \lambda, t\geq t_0$.

Then the fundamental function and the fundamental system of  (\ref{1}) are positive
for $t\geq t_0$.
\end{guess}
{\bf Proof.} It is sufficient to prove the theorem for the case $b(t)\geq 0$.

1) The function $u(t)=a(t)$  is a solution of inequality (\ref{28}).

2) The function $u(t)=a(t)-\lambda b(t)$  is a solution of inequality (\ref{28}).

3) The function $u(t)=\lambda$  is a solution of inequality (\ref{28}).
\qed

\begin{corollary}
Suppose $a(t)\geq 0$ and at least one of the following conditions holds
\vspace{2mm}

1) $a(t)\geq a>0, ~\int_0^{\infty}b^+(s)ds<\infty$,
\vspace{2mm}

2) $a^2(t)\geq 4B$, where 
${\displaystyle B=\limsup_{t\rightarrow\infty}b^+(t)}$,
\vspace{2mm}

3) ${\displaystyle 
\inf_{\lambda>0}\limsup_{t\rightarrow\infty}\frac{1}{\lambda}
\int_0^t e^{-\int_s^t(a(\xi)-\lambda)d\xi}b^+(s)ds<1}$.

Then there exists $t_0\geq 0$ such that the fundamental function 
and the fundamental system of  (\ref{1}) are positive for $t\geq t_0$.
\end{corollary}
{\bf Proof.} The proof of 1) and 3) is evident.
To prove 2) we assume $\lambda=1/\sqrt{B}$. 
\qed

Now let us compare solutions of the same equation with different right 
hand sides and initial conditions, as well as  of different 
equations. 

\begin{guess}
Suppose $a(t)\geq 0, b(t)\geq 0$, the fundamental function of (\ref{22}) is positive.
Denote by $x(t)$ the solution of the problem 
$$
\ddot{x}(t)+a(t)\dot{x}(t)+b(t)x(t)=f(t),~t\geq t_0,
$$$$
x(t_0)=x_0,~ x'(t_0)=x^{'}_0,
$$
by $v(t)$  the solution of the above problem , where $f(t), x_0,x_0^{'}$
are replaced with $f_1(t), v_0,v_0^{'}$. If $f(t)\geq f_1(t), x_0\geq 
v_0,x_0^{'}\geq v_0^{'}$,
then $x(t)\geq v(t)$. 
\end{guess}
{\bf Proof.}
By Corollary 1 the fundamental system and the fundamental function 
of (\ref{1}) are positive.  
The statement of the theorem follows from solution representation 
formula (\ref{4}).
\qed

\begin{guess}
Suppose $a_1(t)\geq a(t)\geq 0, b(t)\geq b_1(t)\geq 0$, 
the fundamental function of (\ref{22}) is positive.
Denote by $x_1(t),x_2(t), X(t,s)$ the fundamental system and the fundamental 
function of (\ref{1});
by $v_1(t),v_2(t), V(t,s)$ the fundamental system and the fundamental
function of (\ref{35}).
 
Then $x_1(t)\leq v_1(t), x_2(t)\leq v_2(t), X(t,s)\leq V(t,s)$.
\end{guess}
{\bf Proof.}
Denote by $Y(t,s)$ the fundamental function of (\ref{8}), by $Y_1(t,s)$ 
the fundamental function
of (\ref{22}) where $a(t), b(t)$ is replaced by $a_1(t), b_1(t)$, 
respectively. 
By Corollary 1 it is sufficient to prove that $Y_1(t,s)\geq Y(t,s)$.

The function $Y_1(t,0)$ is the solution of the equation
\begin{equation}
\label{41}
\dot{y}(t)+\int_{0}^t e^{-\int_s^t a_1(\xi)d \xi}b_1(s)y(s)ds=0, ~ y(0)=1,
\end{equation}    
which can be rewritten equation in the form
\begin{equation}
\label{42}
\dot{y}(t)+\int_{0}^t e^{-\int_s^t a_1(\xi)d \xi}b(s)y(s)ds=
\int_{0}^t e^{-\int_s^t a(\xi)d \xi}(b(s)-b_1(s))y(s)ds
\end{equation}
$$
+\int_{0}^t \left(e^{-\int_s^t a(\xi)d \xi}-e^{-\int_s^t a_1(\xi)d \xi}\right)b_1(s)y(s)ds, ~ y(0)=1.
$$
By formula (\ref{26}) for the solution of (\ref{42}) we have
$$
Y_1(t,0)=Y(t,0)+\int_0^t Y(t,s)\int_{0}^s e^{-\int_{\tau}^s a(\xi)d \xi}(b(\tau)-b_1(\tau))y(\tau)d\tau ds
$$$$
+\int_0^t Y(t,s)\int_{0}^s \left(e^{-\int_{\tau}^s a(\xi)d \xi}-e^{-\int_{\tau}^s a_1(\xi)d \xi}\right)b_1(\tau)y(\tau)d\tau ds.
$$
Since $a(t)-a_1(t)\geq 0,~ b(t)-b_1(t)\geq 0$, then $Y_1(t,0)\geq Y(t,0)$. 
The inequality $Y_1(t,s)\geq Y(t,s)$ is obtained similarly. 
\qed

\begin{corollary}
Suppose
$$
a(t)\equiv a>0, ~~ b(t)\geq b>0, ~~ a^2-4b\geq 0.
$$
Then the fundamental function $X(t,s)$ of (\ref{1}) is positive and for this function
we have the following estimation
$$
\int_0^t X(t,s)ds\leq\frac{1}{b}.
$$
\end{corollary}
{\bf Proof.}
By Theorem 5 we have $0<X(t,s)\leq V(t,s)$, where $V(t,s)$ is the 
fundamental function of equation (\ref{18}).
Application of Lemma 2 (parts 1 and 3) completes the proof.
\qed

\section{Stability}
\setcounter{equation}{0}

\begin{guess}
Suppose $a(t) \geq \alpha>0, b(t)\geq \beta>0$. If the fundamental function of 
(\ref{1})
 is positive then equation (\ref{1}) is exponentially stable.
\end{guess}
{\bf Proof.}
Consider (\ref{1}) with initial conditions $x(0)=0, \dot{x}(0)=0$.
Suppose $f\in L_{\infty}[0,\infty)$. For the solution of this problem we 
have
$x(t)=\int_0^t X(t,s)f(s)ds$.
By Corollary 4
$$
|x(t)|\leq \int_0^t X(t,s)|f(s)|ds=\int_0^t X(t,s)b(s)\frac{|f(s)|}{b(s)}ds
\leq \frac{\|f\|}{\beta} \int_0^t X(t,s)b(s)ds\leq 
\frac{\|f\|}{\beta},
$$
where $\| \cdot \|$ is the sup-norm in $C[0,\infty)$.
Then $x\in C[0,\infty)$.
We have \\
${\displaystyle
\dot{x}(t)=\int_{0}^t e^{-\int_s^t a(\xi)d \xi}b(s)x(s)ds
+\int_{0}^t e^{-\int_s^t a(\xi)d \xi}b(s)f(s)ds}$,
hence $\dot{x}\in C[0,\infty)$.
By Lemma 3  equation (\ref{1}) is 
exponentially stable. 
\qed

\begin{corollary}
Suppose $a(t)\geq \alpha >0, b(t)\geq \beta>0$ and at least one of 
the conditions of Corollary 7 holds.
Then equation (\ref{1}) is exponentially stable.
\end{corollary}

Denote 
$$
\alpha=\liminf_{t\rightarrow\infty}a(t), \beta=\liminf_{t\rightarrow\infty}b(t),
B=\limsup_{t\rightarrow\infty}b(t).
$$

\begin{guess}
Suppose ${\displaystyle 0<\beta\leq B<\frac{1}{2}{\alpha}^2}$. 
Then equation (\ref{1}) is exponentially stable.
\end{guess}
{\bf Proof.}
Without loss of generality we can assume that $b(t)\geq \beta$ for any 
$t\geq 0$ and there exists $\epsilon>0$ such that
\begin{equation}
\label{43}
\frac{\epsilon}{4} {\alpha}^2<b(t)<\left(\frac{1}{2}-
\frac{\epsilon}{4}\right) {\alpha}^2.
\end{equation}
Consider (\ref{1}) with initial conditions $x(0)=0, \dot{x}(0)=0$.
Suppose $f\in L_{\infty}[0,\infty)$ and let us prove that  $x\in 
C[0,\infty)$.

Denote by $X_0(t,s)$ the fundamental function of the equation
\begin{equation}
\label{44}
\ddot{x}(t)+a(t)\dot{x}(t)+\frac{1}{4}{\alpha}^2x(t)=0.
\end{equation}
By Corollary 7(2), $X_0(t,s)>0, t\geq s\geq 0$. 
Equation (\ref{1}) can be rewritten in the form
in the form
\begin{equation}
\label{45}
\ddot{x}(t)+a(t)\dot{x}(t)+\frac{1}{4}{\alpha}^2x(t)
+ \left( b(t)-\frac{1}{4}{\alpha}^2 \right) x(t)=f(t),
\end{equation}
with the initial conditions $x(0)=0, \dot{x}(0)=0$.
Equation (\ref{45}) is equivalent to 
\begin{equation}
\label{46}
x(t)+\int_0^t X_0(t,s) \left( b(s)-\frac{1}{4}{\alpha}^2 
\right) x(s)ds=g(t),
\end{equation}
where
$g(t)=\int_0^t X_0(t,s) f(s)ds$. 
By Theorem 6  equation (\ref{44}) 
is exponentially stable hence $X_0(t,s)$ has an exponential
estimation. Then $g\in L_{\infty}[0,\infty)$.

Equation (\ref{46}) has the form $x+Px=g$. Corollary 2 implies
\begin{eqnarray*}
|(Px)(t)| & \leq & 
\int_0^t X_0(t,s) \left| b(s)-\frac{1}{4}{\alpha}^2 \right| |x(s)|~ds
\\
& \leq & \int_0^t X_0(t,s) \frac{1}{4}{\alpha}^2 \frac
{|b(s)-\frac{1}{4}{\alpha}^2|}{\frac{1}{4}{\alpha}^2}ds \|x\| 
\\
& \leq & \frac{\sup_{t\geq 
0}|b(t)-\frac{1}{4}{\alpha}^2|}{\frac{1}{4}{\alpha}^2}
\|x\|.
\end{eqnarray*}
Inequality
${\displaystyle
\frac{\sup_{t\geq 0}|b(t)-\frac{1}{4}{\alpha}^2|}{\frac{1}{4}{\alpha}^2} \leq
1-\epsilon
}$
is equivalent to (\ref{43}). Hence $\|P\|<1$,  where the norm is
in $C[0,\infty)$. Then $x\in C[0,\infty)$ for a solution of 
(\ref{46}) and therefore for a solution of (\ref{1}). Similar to the proof of Theorem 6 we have $\dot{x}\in
C[0,\infty)$. 
By Lemma 3 equation (\ref{1}) is exponentially stable.
\qed

\noindent
{\bf Example 1.} Consider the following equation
\begin{equation}
\label{47}
\ddot{x}(t)+a\dot{x}(t)+(1+0.99\sin{t})x(t)=0.
 \end{equation}
If  $a> \sqrt{3.98}$ 
then by Theorem 7 equation (\ref{47})  is exponentially stable.
\vspace{2mm}

In all previous results we obtain stability conditions for 
equations with a positive fundamental function and for "small" 
perturbations of such equations. Below we will give stability conditions 
for equations without any positiveness assumptions.

\begin{guess}
Suppose $a(t)\geq \alpha>0, b(t)\geq 0$ and there exist $A>0, B>0$
such that the following condition holds
\begin{equation}
\label{48}
||a(t)-A||\left\|\frac{b}{a}\right\|+||b(t)-B||<
\left\{\begin{array}{ll}
B ~~ ,&A^2\geq 4B,\\
\frac{A\sqrt{4B-A^2}}{4},&A^2< 4B,
\end{array}\right.
\end{equation}
where $\| \cdot \|$ is the norm in the space $L_{\infty}[t_0,\infty)$
for some $t_0\geq 0$. Then equation  (\ref{1}) is exponentially stable.
\end{guess} 
{\bf Proof.}
Without loss of generality we can assume $t_0=0$.
Consider equation (\ref{1}) with the initial conditions $x(0)=0, 
\dot{x}(0)=0$.
Suppose $f\in L_{\infty}[0,\infty)$. Let us prove that  $x, \dot{x}\in 
C[0,\infty)$.
To this end, rewrite equation (\ref{1}) in the form
\begin{equation}
\label{49}
\ddot{x}(t)+A\dot{x}(t)+Bx(t)+(a(t)-A)\dot{x}(t)
+(b(t)-B)x(t)=f(t), ~x(0)=\dot{x}(0)=0.
\end{equation}
Hence 
\begin{equation}
\label{50}
x(t)+\int_0^t Y(t,s) \left[ (a(s)-A)\dot{x}(s)+(b(s)-B)x(s) 
\right]~ds=f_1(t),
\end{equation}
where $Y(t,s)$ is the fundamental function of 
\begin{equation}
\label{51}
\ddot{x}(t)+A\dot{x}(t)+Bx(t)=0,
\end{equation}
$f_1(t)=\int_0^t Y(t,s) f(s) ds.$
Equation (\ref{51}) is exponentially stable
thus $f_1\in L_{\infty}[0,\infty)$.

Equation (\ref{1}) can be rewritten in a different form
\begin{equation}
\label{52}
\dot{x}(t) +\int_0^t e^{-\int_s^t a(\xi)d \xi}b(s)x(s)ds=r(t), ~ t\geq 0,
\end{equation}
where $r(t)=\int_0^t e^{-\int_s^t a(\xi)d \xi}f(s)ds$.
Since $a(t)\geq \alpha>0$, then $r\in L_{\infty}[0,\infty)$.
Substituting $\dot{x}(t)$ from (\ref{52}) to (\ref{50}) we have
\begin{equation}
\label{53}
x(t)-\int_0^t Y(t,s) \left[(a(s)-A)\int_0^s e^{-\int_{\tau}^s
a(\xi)d \xi}b(\tau)x(\tau)d\tau-(b(s)-B)x(s)\right]ds=f_2(t),
\end{equation}
where
$f_2(t)=f_1(t)+\int_0^t Y(t,s) r(s)ds$. Evidently  
$f_2\in L_{\infty}[0,\infty)$.

Equation (\ref{53}) has the form $x-Hx=f_2$, where
\begin{eqnarray*}
\|(Hx)(t)\| & \leq & \int_0^t Y(t,s) \left[ |a(s)-A|\int_0^s 
e^{-\int_{\tau}^s
a(\xi)d \xi}a(\tau)\frac{b(\tau)}{a(\tau)}d\tau+| b(s)-B| \right]ds ~\|x\|
\\
& \leq & K \left( \|a(t)-A \|\left\|\frac{b}{a}\right\|+ \|b(t)-B\| 
\right)\|x\|,
\end{eqnarray*}
where by Lemma 2 the constant $K$ is 
$$
K=\sup_{t\geq 0}\int_0^t 
|Y(t,s)|ds=\left\{\begin{array}{ll}
\frac{1}{B},& A^2\geq 4B,\\
\frac{4}{A\sqrt{4B-A^2}}, & A^2<4B.\\
\end{array}\right.
$$ 
Then $\|H \|<1$, hence $x\in C[0,\infty)$ for the solution $x$ of 
(\ref{53}) and therefore of (\ref{1}). 
As in the proof of Theorem 7, $\dot{x}\in C[0,\infty)$.
By Lemma 3 equation (\ref{1}) is exponentially stable.
\qed

\begin{corollary}
Suppose there exist $a>0, B>0$ such that
\begin{equation}
\label{54}
B-\frac{a\sqrt{4B-a^2}}{4}<m\leq M<B+\frac{a\sqrt{4B-a^2}}{4},
\end{equation}
where ${\displaystyle m=\liminf_{t\rightarrow\infty}b(t)>0,~
M=\limsup_{t\rightarrow\infty}b(t)}$. Then the equation
\begin{equation}
\label{55}
\ddot{x}(t)+a\dot{x}(t)+b(t)x(t)=0
\end{equation}
is exponentially stable.
\end{corollary}
{\bf Proof}  follows from the second inequality (\ref{38}) if we let $A=a$.
\qed

Application of the first inequality (\ref{38}) 
to equation (\ref{55}) gives the same stability conditions
which were obtained by application of Theorem 7.

{\bf Example 2.}
Consider the equation
\begin{equation}
\label{56}
\ddot{x}(t)+(10+\alpha(t)\sin t)\dot{x}(t)+(26+\beta(t)\cos t)x(t)=0,
\end{equation}
where $\alpha(t), \beta(t)$ are measurable functions, such that 
$|\alpha(t)|\leq 1, |\beta(t)|\leq 1$.

If we take $A=10, B=26$, then all condition of Theorem 8 hold.
Then equation (\ref{56}) is exponentially stable.
\vspace{2mm}

We will obtain new stability conditions using the derivative of the 
fundamental function of comparison equations.

\begin{guess}
Suppose $a(t)\geq \alpha>0, b(t)\geq 0$ and there exist $a>0, 
b>0$ such that at least one of the following conditions holds:

1)
${\displaystyle
a^2>4b, ~\|a(t)-a\|\frac{2a}{\sqrt{a^2-4b}(a-\sqrt{a^2-4b})}+\|b(t)-b\|\frac{1}{b}<1.
}$

2)
${\displaystyle
a^2<4b, ~\|a(t)-a\|\frac{2(a+\sqrt{4b-a^2})}{a\sqrt{4b-a^2}}+\|b(t)-b\|\frac{4}{a\sqrt{4b-a^2}}<1.
}$

3)
${\displaystyle
a^2=4b, ~\|a(t)-a\|\frac{2}{\sqrt{b}}+\|b(t)-b\|\frac{1}{b}<1,
}$
\\
where $\| \cdot \|$ is the norm in the space $L_{\infty}[t_0,\infty)$
for some $t_0\geq 0$. 
Then equation  (\ref{1}) is exponentially stable.
\end{guess} 
{\bf Proof.}
Without loss of generality we can assume $t_0=0$.
Suppose $x(t)$ is a solution of (\ref{1})
with initial conditions $x(0)=x'(0)=0$. Denote $z(t)=\ddot{x}(t)+a\dot{x}(t)+bx(t), ~Y(t,s)$
is the fundamental function of the equation $\ddot{x}(t)+a\dot{x}(t)+bx(t)=0$. Then
\begin{equation}
\label{57}
x(t)=\int_0^t Y(t,s)z(s)ds,~ \dot{x}(t)=\int_0^t Y_t^{'}(t,s)z(s)ds.
\end{equation}
Equation (\ref{1}) is equivalent to the equation
\begin{equation}
\label{58}
\ddot{x}(t)+a\dot{x}(t)+b(t)x(t)+(a(t)-a)\dot{x}(t)+(b(t)-b)x(t)=f(t).
\end{equation}
After substituting $z(t)$ and (\ref{57}) into (\ref{58})
we have the following equation 
\begin{equation}
\label{59}
z(t)+(a(t)-a)\int_0^t Y_t^{'}(t,s)z(s)ds+(b(t)-b)\int_0^t Y(t,s)z(s)ds=f(t).
\end{equation}
Equation (\ref{59}) has the form $z+Hz=f$. For the norm of the operator 
$H$ in the space 
$L_{\infty}[0,\infty)$ we have
$$
\|H\|\leq \|a(t)-a\|\left\|\int_0^t 
|Y_t^{'}(t,s)|ds\right\|+\|b(t)-b\|\left\|\int_0^t |Y(t,s)|ds\right\|.
$$
By Lemma 2 conditions 1)-3) of the theorem imply $\|H\|<1$. Hence for 
the solution of (\ref{59})
we have $z\in L_{\infty}[0,\infty)$. Equation $\ddot{x}(t)+a\dot{x}(t)+bx(t)=0$ is exponentially stable.
Equalities (\ref{57}) imply that the solution of (\ref{1}) and its derivative are bounded functions.
Then by Lemma 3 equation (\ref{1}) is exponentially stable.
\qed 

\noindent
{\bf Example 3.} Consider the following equation
\begin{equation}
\label{60}
\ddot{x}(t)+\dot{x}(t)+(b+\sin{t})x(t)=0.
\end{equation}
If $b>4.25$ then by Theorem 9(2) ($a=1$, $b$ is the same as in the 
equation)
equation (\ref{60})  is exponentially stable.
Theorem 8 gives the same result, if we take $a=1, B=b$. 
Theorem 7 does not give any stability condition for this equation.
\vspace{2mm}

\section{Zones of Lyapunov's Stability}

In the following, let us assume that $a(t)$
is an absolutely continuous function.

\setcounter{equation}{0}

It is known that the substitution 
\begin{equation}
\label{61}
x(t)=z(t)\exp \left\{ -\frac{1}{2}\int_{0}^{t}a(s)ds\right\} 
\end{equation}
transforms the homogeneous equation 
\begin{equation}
\label{62}
\pounds x(t)\equiv \ddot{x}(t)+a(t)\dot{x}(t)+b(t)x(t)=0,\;t\in [ 0,\infty ),
\end{equation}
into the equation 
\begin{equation}
\label{63}
\ddot{z}(t)+p(t)z(t)=0,\;t\in [ 0,\infty ),
\end{equation}
where 
\begin{equation}
\label{64}
p(t)=b(t)-\frac{a^{2}(t)}{4}-\frac{a^{\prime }(t)}{2}.
\end{equation}
Obviously coefficient $p(t)$ can be presented as $
p(t)=p^{+}(t)-p^{-}(t),$ where $p^{+}(t)\geq 0$ and $p^{-}(t)\geq 0.$

Consider now equation (\ref{62}) with an $\omega$-periodic
coefficients $a(t)$ and $b(t).$

It is known from the works of the well known mathematicians  
Zhukovskii \cite{zhukovskii}, Kre\u{\i}n \cite{krein} and Yakubovich 
\cite{yakubovich1} that there is a deep connection between the problem of 
the Lyapunov's stability and the nonoscillation intervals. We propose the 
following statement.

\begin{guess} Assume that $a(t+\omega )=a(t),\;b(t+\omega )=b(t)\;$%
for $t\in [0,\infty )$ and 
\begin{equation}
\label{65}
\int_{0}^{\omega }p(t)dt>0,
\end{equation}
where $p(t)$ is defined in (\ref{64}),
and at least one of the following three conditions holds:

1) $a(t)\geq \int_{0}^{t}b^{+}(s)ds$  for $t\in [0,\omega ];$

2) $a(t)\geq 0$ and there exists $\lambda >0$ such that $
\int_{0}^{t}\exp \left\{ -\int_{s}^{t}(a(\xi )-\lambda )d\xi 
\right\}
b^{+}(s)ds\leq \lambda $  for $t\in [0,\omega ];$

3) $\int_{0}^{\omega }p^{+}(t)dt\leq \frac{4}{\omega }.$

Then all solutions of homogeneous equation (\ref{62}) tend to
zero when $t\to \infty$, if $\int_{0}^{\omega }a(t)dt>0$,
and all solutions are bounded if $\int_{0}^{\omega }a(t)dt=0$.
\end{guess}
{\bf Proof.} It is known \cite{levin} that if $[0,\omega ]$ is a 
nonoscillation interval for (\ref{62}),
where $\omega $ is the period of the coefficient $p(t)$, then
condition (\ref{65}) garantees that all solutions of equation (\ref{62}) 
are bounded. Each of the conditions 1)-3) yields that $[0,\omega ]$ is a 
nonoscillation interval.
The conditions on the integral of the function $a(t)$ and reference to
the substitution (\ref{61}) completes the proof.
\qed

\section{Floquet Theory and Stability}
\setcounter{equation}{0}

Consider now the equation 
\begin{equation}
\label{66}
\pounds x(t)\equiv \ddot{x}(t)+a(t)\dot{x}(t)+b(t)x(t)=0,\;t\in [0,\infty ),
\end{equation}
with $\omega$-periodic coefficients $a(t+\omega 
)=a(t)$, $b(t+\omega)=b(t)$. For this equation there exist solutions 
satisfying the condition 
\begin{equation}
\label{67}
x(t+\omega )=\lambda x(t).
\end{equation}

The foundations and applications of the Floquet theory were presented in the
book by Yakubovich and Starzhinski\u{\i} \cite{yakubovich}.
Using the Floquet theory for ordinary differential equations write the
equation for $\lambda $: 
\begin{equation}
\label{68}
\lambda ^{2}-(x_{1}(\omega )+x_{2}^{\prime }(\omega ))\lambda +W(\omega )=0,%
\end{equation}
where $x_{1}$ and $x_{2}$ are two solutions of the equation 
(\ref{66})  such that $x_{1}(0)=1$, $x_{1}^{\prime }(0)=0$, $x_{2}(0)=0$, 
$x_{2}^{\prime}(0)=1$. Denote by 
$$
W(t)=\det \left( 
\begin{array}{cc}
x_{1}(t) & x_{2}(t) \\ 
x_{1}^{\prime }(t) & x_{2}^{\prime }(t)
\end{array}
\right) 
$$
the Wronskian of the fundamental system of (\ref{66}). Obviously 
$W(0)=1$.

If $\lambda _{1}$ is a real root of equation (\ref{67}), then the 
corresponding solution 
of equation (\ref{66}) 
has the representation 
\begin{equation}
\label{69}
y(t)=g(t)\exp \left( \frac{\ln |\lambda _{1}|}{\omega }t\right) ,
\end{equation}
where $g$ is $\omega$-periodic if $\lambda _{1}>0$ and is $
2\omega $-periodic function if $\lambda _{1}<0$. If equation (\ref{68}) has two
complex roots $\lambda _{1}=|\lambda _{1}|\exp (i\theta )$ and $\lambda
_{2}=|\lambda _{1}|\exp (-i\theta )$, then the corresponding solutions 
of equation (\ref{66}) have the form
\begin{equation}
\label{70}
y_{1}(t)=\left( g_{1}(t)\cos \frac{\theta t}{\omega }-g_{2}(t)\sin \frac{%
\theta t}{\omega }\right) \exp \left( \frac{\ln |\lambda _{1}|}{\omega }%
t\right) ,
\end{equation}
\begin{equation}
\label{71}
y_{2}(t)=\left( g_{2}(t)\cos \frac{\theta t}{\omega }+g_{1}(t)\sin \frac{%
\theta t}{\omega }\right) \exp \left( \frac{\ln |\lambda _{1}|}{\omega }%
t\right),
\end{equation}
where $g_{1}$ and $g_{2}$ are $\omega $-periodic functions.

\begin{guess} Assume that equation (\ref{66}) is oscillatory and the
distance between zeros of its solutions is different from $2\omega$.
Then the following statements are valid.

a) Equation (\ref{66}) is exponentially stable if  
\begin{equation}
\label{72}
\int_{0}^{\omega }a(t)dt>0.
\end{equation}

b) The fundamental solutions of equation (\ref{66}) are of the form
(\ref{69}), where $\left| \lambda _{1}\right| >1$  if 
\begin{equation}
\label{73}
\int_{0}^{\omega }a(t)dt<0.
\end{equation}

c) If  
\begin{equation}
\label{74}
\int_{0}^{\omega }a(t)dt=0
\end{equation}
then the fundamental solutions of equation (\ref{66}) are bounded. 
\end{guess}
{\bf Proof.}
It follows from the classical formula of Ostrogradskii that condition (\ref{72})
implies the inequality $W(\omega )<1,$ condition (\ref{73}) implies the
inequality $W(\omega )>1,$ and condition (\ref{74}) implies that $W(\omega )=1.$
The condition that the distance between zeros of solutions of (\ref{66}) is
different from $2\omega $ excludes the existence of real roots of equation
(\ref{68}). In this case the inequality $W(\omega )<1$ implies that $\left|
\lambda _{1}\right| <1,$ the equality $W(\omega )=1$ implies that $\left|
\lambda _{1}\right| =1,$ and the inequality $W(\omega )>1$ implies that 
$\left| \lambda _{1}\right| >1.$ Now the representation of solutions 
(\ref{70}),(\ref{71}) completes the 
proof.
\qed

\noindent
{\bf Remark.} The condition that the distance between zeros
are different from $2\omega $\ is essential as the following example
demonstrates. 

\noindent{\bf Example 4.}
Consider the equation 
\begin{equation}
\label{75}
\pounds x(t)\equiv x^{\prime \prime }(t)+\frac{2\sin ^{2}t+\cos t\sin t}
{1+\cos t\sin t}x^{\prime }(t)+\frac{\sin ^{2}t-\cos t\sin t}{1+\cos t\sin 
t}x(t)=0,\;t\in [ 0,\infty ).
\end{equation}
Inequality (\ref{72}) for the coefficient 
$a(t)=\frac{2\sin ^{2}t + \cos t\sin t}
{1+\cos t\sin t}$\ is fulfilled with $\omega=\pi$, but this equation is 
not exponentially stable: its fundamental system is $x_{1}=e^{-t}\cos t$ 
and $x_{2}=\sin t$.

\bigskip

Using the substitution (\ref{61}), we again obtain the equation 
\begin{equation}
\label{76}
z^{\prime \prime }(t)+p(t)z(t)=0,~~t\in [0,\infty ),
\end{equation}
where 
\begin{equation}
\label{77}
p(t)=b(t)-\frac{a^{2}(t)}{4}-\frac{a^{\prime }(t)}{2},
\end{equation}
Evidently zeros of the solution $x$ of equation (\ref{66}) and the
corresponding solution $z$ of the equation (\ref{63}) coinside. Let us denote 
\begin{equation}
\label{78}
P=ess\inf_{t\in [0,\omega ]}p(t),\;Q=ess\sup_{t\in [0,\omega
]}p(t).
\end{equation}
Estimating distances  between two adjacent zeros (i.e. nonoscillation
intervals) from below and from above we get the following result.

\begin{guess} Suppose $P>0,$ there exists a positive integer $k$ such that
$\frac{k-1}{k}<\sqrt{\frac{P}{Q}}$ and   
\begin{equation}
\label{79}
\omega \in \left( 0,\frac{\pi }{2\sqrt{Q}}\right] \cup \left( 
\frac{1}{2}\frac{\pi }{%
\sqrt{P}},\frac{\pi }{\sqrt{Q}}\right) \cup ...\cup \left( \frac{k-1}{2}%
\frac{\pi }{\sqrt{P}},\frac{k\pi }{2\sqrt{Q}}\right) .
\end{equation}
Then equation (\ref{66}) is oscillatory and distance between zeros of its
solutions is different from $2\omega .$
\end{guess}
{\bf Proof. }If equation (\ref{68}) has real roots, then there exist such  
zeros $t_{0},t_{1}$ of a solution $x(t)$ that the distance between 
$t_{0}$ and $t_{1}$ equals $\omega $ or $2\omega$. We will reject this 
possibility, since the distance between zeros of $g(t)$ in (\ref{69})
cannot be $2\omega$.

Assume that $x(t_{0})=0.$ We use the functions $v=\sin \sqrt{Q}%
(t-t_{0})$ in the first assertion of Theorem A to get that the
spectral radius of the operator $K_{t_{0},t_{0}+\frac{\pi }{2\sqrt{Q}}%
}:C_{[t_{0},t_{0}+\frac{\pi }{2\sqrt{Q}}]}\rightarrow C_{[t_{0},t_{0}+\frac{%
\pi }{2\sqrt{Q}}]}$ defined by the equality 
\begin{equation}
\label{80}
K_{t_{0},t_{0}+\frac{\pi }{2\sqrt{Q}}}x(t)=-\int_{t_{0}}^{t_{0}+\frac{\pi }{2%
\sqrt{Q}}}G_{t_{0},t_{0}+\frac{\pi }{2\sqrt{Q}}}(t,s)p(s)x(s)ds.
\end{equation}
where $G_{t_{0,}t_{0}+\frac{\pi }{2\sqrt{Q}}}(t,s)$\ is the Green's function
of the problem 
\begin{equation}
\label{81}
x^{\prime \prime }(t)+a(t)x^{\prime }(t)+b(t)x(t)=f(t),\;t\in 
[0,\omega ],
x(t_{0})=0,\;x\left( t_{0}+\frac{\pi }{2\sqrt{Q}}\right)
=0,
\end{equation}
is less than one. 

Applying Theorem 5.4 in \cite{krasnoselskii}, p. 81, we obtain that the
spectral radius of the operator
$
K_{t_{0},t_{0}+\frac{\pi }{2\sqrt{P}}}:C_{[t_{0},t_{0}+\frac{\pi }{2\sqrt{P}}
]}\rightarrow C_{[t_{0},t_{0}+\frac{\pi }{2\sqrt{P}}]},$ defined by the
equality 
\begin{equation}
\label{82}
K_{t_{0},t_{0}+\frac{\pi }{2\sqrt{P}}}x(t)=-\int_{t_{0}}^{t_{0}+\frac{\pi }{2
\sqrt{P}}}G_{t_{0},t_{0}+\frac{\pi }{2\sqrt{P}}}(t,s)p(s)x(s)ds.
\end{equation}
where $G_{t_{0},t_{0}+\frac{\pi }{2\sqrt{P}}}(t,s)$ is the Green's function
of the problem 
\begin{equation}
\label{83}
x^{\prime \prime }(t)+a(t)x^{\prime }(t)+b(t)x(t)=f(t),\;t\in [0,\omega ],
x(t_{0})=0,\;x\left( t_{0}+\frac{\pi }{2\sqrt{P}}\right)
=0,
\end{equation}
is greater or equal to one. Moving the point $t_{0}$ we obtain that there
are no zeros in the zones defined by (\ref{79}).
\qed

{\bf Acknowlegment}: the research of the first and the third author were 
supported by The Israel Science Foundation (grant No. 828/07), the second
author was partially supported by National Sciences and Engineering 
Research Council of Canada.

\end{document}